\documentclass[11pt]{amsart}
\author{Paul Pollack}
\address{Department of Mathematics\\ University of Georgia\\ Athens, GA 30602}
\email{pollack@uga.edu}
\thanks{The author is supported by the National Science Foundation (NSF) under award DMS-2001581.}

\title{Numbers which are orders only of cyclic groups}
\usepackage{amsmath,amssymb,amsthm,microtype,geometry,mathscinet,parskip,enumitem}
\setlist[itemize]{parsep=0pt}
\setlist[enumerate]{parsep=0pt}
\usepackage[utf8]{inputenc}

\DeclareMathAlphabet{\curly}{U}{rsfs}{m}{n}
\newtheorem{thm}{Theorem}[section]

\newtheorem{lem}[thm]{Lemma}

\theoremstyle{remark}

\begin{document}
\renewcommand{\labelenumi}{(\roman{enumi})}
\def\Ll{\mathcal{L}}
\def\N{\mathbb{N}}
\def\Aa{\mathcal{A}}
\def\Q{\mathbb{Q}}
\newcommand\rad{\mathrm{rad}}
\def\Z{\mathbb{Z}}
\def\F{\mathbb{F}}
\def\R{\mathbb{R}}
\def\C{\mathbb{C}}
\def\Pp{\mathcal{P}}
\def\Ss{\mathcal{S}}
\newcommand\Li{\mathrm{Li}}
\begin{abstract} We call $n$ a \textsf{cyclic number} if every group of order $n$ is cyclic. It is implicit in work of Dickson, and explicit in work of Szele, that $n$ is cyclic precisely when $\gcd(n,\phi(n))=1$. With $C(x)$ denoting the count of cyclic $n\le x$, Erd\H{o}s proved that
\[ C(x) \sim e^{-\gamma} x/\log\log\log{x}, \quad\text{as $x\to\infty$}. \]
We show that $C(x)$ has an asymptotic series expansion, in the sense of Poincar\'e, in descending powers of $\log\log\log{x}$, namely
\[ \frac{e^{-\gamma} x}{\log\log\log{x}} \left(1-\frac{\gamma}{\log\log\log{x}} + \frac{\gamma^2 + \frac{1}{12}\pi^2}{(\log\log\log{x})^2} - \frac{\gamma^3 +\frac{1}{4} \gamma \pi^2 + \frac{2}{3}\zeta(3)}{(\log\log\log{x})^3} + \dots \right). \]
\end{abstract}

\maketitle

\section{Introduction}
Call the positive integer $n$ \emph{cyclic} if the cyclic group of order $n$ is the unique group of order $n$. For instance, all primes are cyclic numbers. It is implicit in work of Dickson \cite{dickson05}, and explicit in work of Szele \cite{szele47}, that $n$ is cyclic precisely when $\gcd(n,\phi(n))=1$, where $\phi(n)$ is Euler's totient. (In fact, this criterion had been stated as ``evident'' already by Miller in 1899 \cite[p.\ 235]{miller99}.) If $C(x)$ denotes the count of cyclic numbers $n\le x$, Erd\H{o}s proved in \cite{erdos48} that
\begin{equation}\label{eq:erdosasymptotic} C(x) \sim e^{-\gamma} x/\log\log\log{x}, \end{equation}
as $x\to\infty$, where $\gamma$ is the Euler--Mascheroni constant. Thus, the relative frequency of cyclic numbers decays to $0$ but ``with great dignity'' (Shanks).

Several authors have investigated analogues of \eqref{eq:erdosasymptotic} for related counting functions from enumerative group theory. See, for example, \cite{mays79, MM84, warlimont85, srinivasan87, EMM87, EM88, NS88, srinivasan91, NP18}. Our purpose in this note is somewhat different; we aim to refine the formula \eqref{eq:erdosasymptotic}. Begunts \cite{begunts01}, optimizing the method of \cite{erdos48}, showed that $C(x)$ is given by $e^{-\gamma}/\log\log\log{x}$ up to a multiplicative error of size $1+O(\log\log\log\log{x}/\log\log\log{x})$  (the same result appears as Exercise 2 on p.\ 390 of \cite{MV07}). We improve this as follows.

\begin{thm}\label{thm:main} The function $C(x)$ admits an asymptotic series expansion, in the sense of Poincar\'e {\rm(}see \cite[\S1.5]{dB81}{\rm)}, in descending powers of $\log\log\log{x}$. \emph{Precisely:} There is a sequence of real numbers $c_1, c_2, c_3, \dots$ such that, for each fixed positive integer $N$ and all large $x$,
\begin{multline*} C(x) = \frac{e^{-\gamma} x}{\log\log\log{x}} \left(1+\frac{c_1}{\log\log\log{x}} + \frac{c_2}{(\log\log\log{x})^2} + \dots + \frac{c_N}{(\log\log\log{x})^{N}}\right)\\ + O_{N}\left(\frac{x}{(\log\log\log{x})^{N+2}}\right). \end{multline*}
\end{thm}

Our proof of Theorem \ref{thm:main} yields the following explicit determination of the constants $c_k$. Write the Taylor series for the $\Gamma$-function, centered at $1$, in the form $\Gamma(1+z) = 1 + C_1 z + C_2 z^2 + \dots$.
Then the coefficients $c_1, c_2, \dots$ are determined by the formal relation
\[ 1 + c_1 z + c_2 z^2 + c_3 z^3 + \dots = \exp(0! C_1 z + 1! C_2 z^2 + 2! C_3 z^3 + \dots). \]

 For computations of the $C_k$ and $c_k$, it is useful to recall that \begin{equation}\label{eq:digamma} \Gamma(1+z) = \exp\left(-\gamma z + \sum_{k=2}^{\infty}\frac{(-1)^k}{k} \zeta(k) z^k\right). \end{equation}
(This is one version of a well-known expansion for the digamma function; see, e.g., entries 5.7.3 and 5.7.4 in \cite{NIST10}.) The first few $c_k$ are given by
\[ c_1 = -\gamma, \quad c_2 = \gamma^2 + \frac{1}{2}\zeta(2) = \gamma^2 + \frac{\pi^2}{12}, \quad c_3 = -\left(\gamma^3+\frac{1}{4}\gamma \pi^2 + \frac{2}{3}\zeta(3)\right). \]
Owing to \eqref{eq:digamma}, each $c_k$ belongs to the ring $\Q[\gamma,\zeta(2),\zeta(3),\dots\zeta(k)]$. From the fact that the coefficients of $\log \Gamma(1+z)$ are alternating in sign, one deduces that both the $C_k$ and the $c_k$ are alternating as well. Moreover, $$|c_k| \ge (k-1)!|C_k| \ge (k-1)! \zeta(k)/k \ge (k-1)!/k$$ for each $k\ge 2$. It follows that the  series $1+c_1/\log\log\log{x} + c_2/(\log\log\log{x})^2 + \dots$ is purely an asymptotic series, in that it diverges for all values of $x$.

The proof of Theorem \ref{thm:main} has many ingredients in common with the related work cited above (see also \cite{PP20, pollack20}). But we must be more careful about error terms than in earlier papers, and somewhat delicate bookkeeping is required to wind up with a clean result.

\subsection*{Notation} The letters $p$ and $q$ are reserved for primes. We use $K_0, K_1, K_2$, etc.\ for absolute positive constants. To save space, we write $\log_k$ for the $k$th iterate of the natural logarithm.

\section{Lemmata}
We will use Mertens' theorem in the following form, which is a consequence of the prime number theorem with the classical $x\exp(-K_0\sqrt{\log{x}})$ error estimate of de la Vall\'ee Poussin.

\begin{lem}\label{lem:mertens} There is an absolute constant $c$ such that, for all $X\ge 3$,
\[ \sum_{p\le X} \frac{1}{p} = \log_2{X} + c + O(\exp(-K_1\sqrt{\log X})). \]
Moreover, for all $X\ge 3$,
\[ \prod_{p \le X}\left(1-\frac{1}{p}\right) = \frac{e^{-\gamma}}{\log{X}}\left(1+O(\exp(-K_2\sqrt{\log X}))\right). \]
\end{lem}

The following sieve result is a special case of \cite[Theorem 7.2]{HR74}.

\begin{lem}\label{lem:sieve} Suppose that $X \ge Z\ge 3$. Let $\mathcal{P}$ be a set of primes not exceeding $Z$. The number of $n\le X$ coprime to all elements of $\mathcal{P}$ is
\[ X \prod_{p \in \mathcal{P}} \left(1-\frac{1}{p}\right) \left(1 + O\left(\exp\left(-\frac{1}{2}\frac{\log X}{\log Z}\right)\right)\right). \]
\end{lem}

The final estimate of this section was proved independently by Pomerance (see Remark 1 of \cite{Pomerance77}) and Norton (see the Lemma on p.\ 699 of \cite{norton76}).

\begin{lem}\label{lem:reciprocals} For every positive integer $m$ and every $X\ge 3$,
%\[ S(X;m) = \sum_{\substack{p \le X \\ p \,\equiv\,1\kern-3pt\!\!\!\pmod{m}}} \frac{1}{p}. \]
%Then $$S(X;m) = \frac{\log_2{X}}{\phi(m)} + O\left(\frac{\log{(2m)}}{\phi(m)}\right).$$
\[ \sum_{\substack{p \le X \\ p \,\equiv\,1\!\!\!\!\pmod{m}}} \frac{1}{p}= \frac{\log_2{X}}{\phi(m)} + O\left(\frac{\log{(2m)}}{\phi(m)}\right).\]

\end{lem}

\section{Proof of Theorem \ref{thm:main}}
\subsection{Outline}\label{sec:outline} We summarize the strategy of the proof, deferring the more intricate calculations to later sections. Put
\[ y = \frac{\log_2 x}{2\log_3 x} \quad{\text{and}}\quad z = (\log_2 x)\cdot \exp(\sqrt{\log_3 x}). \]

Let us call the prime $p$ a \textsf{standard divisor of $\gcd(n,\phi(n))$} if there is a prime $q \le x^{1/\log_2 x}$ dividing $n$ with $q\equiv 1\pmod{p}$. Clearly, each standard divisor of $(n,\phi(n))$ is a divisor of $\gcd(n,\phi(n))$.

Let $\Ss_0$ be the set of $n\le x$ with no prime factor in $[2,y]$. For each positive integer $k$, let $\Ss_k$ be the set of $n \in \Ss_0$ having exactly $k$ distinct prime factors from the interval $(y,z]$, all of which divide $n$ to the first power only, and at least one of which is a standard divisor of $\gcd(n,\phi(n))$. We will estimate $C(x)$ by
\begin{equation}\label{eq:ourestimate} \#\left(\Ss_0\setminus \bigcup_{1\le k\le \log_3 x} \Ss_k\right) = \#\Ss_0 - \sum_{1 \le k \le \log_3 x} \#\Ss_k. \end{equation}

Suppose $n$ is counted by $C(x)$ but not by \eqref{eq:ourestimate}. Then $n$ has a prime factor $p\le y$. Since $n$ is counted by $C(x)$, it must be that $p\nmid \phi(n)$, so that $n$ is not divisible by any $q\equiv 1\pmod{p}$. By Lemma \ref{lem:sieve}, for a given $p$ the number of those $n\le x$ is $\ll x \prod_{q\le x,\ q \equiv 1\pmod{p}}(1-1/q) \le x \exp(-\sum_{q\le x,\, q\equiv 1\pmod{p}}1/q)$. And by Lemma \ref{lem:reciprocals}, $$ \sum_{\substack{q \le x\\q\equiv 1\!\!\!\!\pmod{p}}} \frac{1}{q} = \frac{1}{p-1} \log_2 x + O(1) \ge 2\log_3 x + O(1).$$ Thus, the number of $n$ corresponding to a given $p$ is $\ll x\exp(-2\log_3 x) = x/(\log_2 x)^2$. Summing on $p\le y$, we deduce that the total number of $n$ counted by $C(x)$ but not \eqref{eq:ourestimate} is $O(x/\log_2 x)$.

Working from the opposite side, suppose that $n$ is counted by \eqref{eq:ourestimate} but not by $C(x)$. Then at least one of the following holds:
\begin{enumerate}
\item there is a prime $p > y$ for which $p^2 \mid n$,
\item there is a prime $p > z$ that divides $n$ and $\phi(n)$,
\item there is a prime $p$ in $(y,z]$ dividing $n$ and a prime $q\equiv 1\pmod{p}$ dividing $n$ with $q > x^{1/\log\log{x}}$,
\item $n$ has more than $\log_3 x$ different prime factors in $(y,z]$.
\end{enumerate}
The number of $n\le x$ for which (i) holds is $\ll x\sum_{p>y} 1/p^2 \ll x/y\log{y} \ll x/\log_2 x$. In order for (ii) to hold but (i) to fail, there must be a prime $q\equiv 1\pmod{p}$ dividing $n$. Clearly, there are most $x/pq$ such $n$ corresponding to a given $p, q$. Thus, the number of $n$ that arise this way is
\[ \ll x\sum_{p > z} \frac{1}{p} \sum_{\substack{q \le x\\q\equiv 1\!\!\!\!\pmod{p}}} \frac{1}{q} \ll x\sum_{p>z} \frac{\log_2 x + \log{p}}{p^2} \ll \frac{x \log_2{x}}{z} = \frac{x}{\exp(\sqrt{\log_3 x})}. \]
For similar reasons, the number of $n\le x$ for which (iii) holds is
\[ \ll x\sum_{y < p \le z} \frac{1}{p} \sum_{\substack{x^{1/\log_2 x} < q \le x \\ q\equiv 1\!\!\!\!\pmod{p}}} \frac{1}{q} \ll x\sum_{p > y} \frac{\log_3 x}{p^2} \ll x\frac{\log_3 x}{\log_2 x}. \]
To handle (iv), observe that $\sum_{y < p \le z}1/p \le K_3/\sqrt{\log_3 x} < 1/2$ for large values of $x$. Thus, the number of $n\le x$ for which (iv) holds is (crudely) at most
\[  x\sum_{k > \log_3 x} \bigg(\sum_{y < p \le z} 1/p\bigg)^k \le 2x (K_3/\sqrt{\log_3 x})^{\log_3 x} \le x/\log_2 x. \]

Collecting estimates, we conclude that
\[ C(x) = \#\left(\Ss_0\setminus \bigcup_{1\le k\le \log_3 x} \Ss_k\right) + O(x/\exp(\sqrt{\log_3 x})). \]
Since the error term is $O_N(x/(\log_3 x)^{N+2})$ for any fixed $N$, for the sake of proving  Theorem \ref{thm:main} we may replace $C(x)$ by $\#(\Ss_0\setminus \bigcup_{1\le k\le\log_3 x} \Ss_k)$.

In \S\ref{sec:sskestimate} we prove suitable estimates for the numbers $\#\Ss_k$ and in \S\ref{sec:DayNewMa} we tie everything together and complete the proof of Theorem \ref{thm:main}.
\subsection{Estimating $\#\Ss_k$}\label{sec:sskestimate} The case $k=0$ is easy to dispense with. By Lemmas \ref{lem:mertens} and \ref{lem:sieve},
\begin{equation}\label{eq:s0estimate} \#\Ss_0 = e^{-\gamma}\frac{x}{\log y} + O(x/\exp(K_4\sqrt{\log_3 x})). \end{equation}

Now suppose that $1 \le k \le \log_3 x$. In order for the integer $n\le x$ to be counted by $\Ss_k$, it is necessary and sufficient than $n=p_1 \cdots p_k m$, where (a) $p_1,\dots,p_k$ are distinct primes belonging to $(y,z]$, (b) the integer $m$ is free of prime factors in $[2,z]$, and (c) $m$ has a prime factor $q \le x^{1/\log_2 x}$ with $q\equiv 1 \pmod{p_i}$ for some $i=1,2,\dots,k$.

Fix distinct primes $p_1,\dots,p_k \in (y,z]$. We will count the number of $n \in \Ss_k$ for which $p_1,\dots,p_k$ are the
prime divisors of $n$ in $(y,z]$. To get at this, we count all $n=p_1\dots p_k m\le x$ where condition (b) holds and then subtract the contribution from $n$ for which (b) holds but (c) fails. By Lemma \ref{lem:sieve}, this is approximately
\begin{equation}\label{eq:skapproximationfirst} \frac{x}{p_1\cdots p_k} \prod_{p \le z}\left(1-\frac{1}{p}\right) \Bigg(1- \prod_{\substack{z < q\le x^{1/\log_2 x} \\ q\equiv 1\!\!\!\!\!\pmod{p_i}\,\text{for some $i$}}}\left(1-\frac{1}{q}\right) \Bigg). \end{equation}
In fact, taking $X=x/p_1 \cdots p_k$ (which exceeds $x^{1/2}$) and $Z = x^{1/\log\log{x}}$ in Lemma \ref{lem:sieve}, we see that the error in this approximation is (very crudely) bounded by $O(x/(p_1\dots p_k \log_2 x))$.

Now we replace $\prod_{p\le z} (1-1/p)$ in \eqref{eq:skapproximationfirst} with $e^{-\gamma}/\log{z}$. This introduces another error of size $x/(p_1\cdots p_k \exp(K_5 \sqrt{\log_3 x}))$.

It remains to estimate the product over $q$ in \eqref{eq:skapproximationfirst}. We have that
\begin{align*} \prod_{\substack{z < q\le x^{1/\log_2 x} \\ q\equiv 1\!\!\!\!\!\pmod{p_i}\,\text{for some $i$}}}\left(1-\frac{1}{q}\right) &= \exp\Bigg(-\sum_{\substack{z < q\le x^{1/\log_2 x} \\ q\equiv 1\!\!\!\!\!\pmod{p_i}\,\text{for some $i$}}}\frac{1}{q} + O\left(\sum_{q > z}\frac{1}{q^2}\right)\Bigg) \\
&= \exp\Bigg(-\sum_{\substack{z < q\le x^{1/\log_2 x} \\ q\equiv 1\!\!\!\!\!\pmod{p_i}\,\text{for some $i$}}}\frac{1}{q}\Bigg)(1+O(1/z)).
\end{align*}
Continuing, we observe that
\[ \sum_{\substack{z < q\le x^{1/\log_2 x} \\ q\equiv 1\!\!\!\!\!\pmod{p_i}\,\text{for some $i$}}}\frac{1}{q}  = \sum_{i=1}^{k} \sum_{\substack{z < q \le x^{1/\log_2 x} \\ q \equiv 1\!\!\!\!\!\pmod{p_i}}} \frac{1}{q} + O\Bigg(\sum_{1 \le i < j \le k}\sum_{\substack{z < q \le x^{1/\log_2 x} \\ q \equiv 1\!\!\!\!\!\pmod{p_i p_j}}}\frac{1}{q}\Bigg), \]
and that the $O$-term here is
\[ \ll \sum_{1 \le i < j \le k}\frac{\log_2 x}{p_i p_j} \ll \binom{k}{2}  \frac{(\log_3 x)^2}{\log_2 x} \ll \frac{(\log_3 x)^4}{\log_2 x}. \]
Moreover,
\begin{align*} \sum_{i=1}^{k} \sum_{\substack{z < q \le x^{1/\log_2 x} \\q \equiv 1\!\!\!\!\pmod{p_i}}} \frac{1}{q} &= \sum_{i=1}^{k} \left(\frac{\log_2 x}{p_i-1} + O\left(\frac{\log_3 x}{p_i}\right)\right)
\\&= \sum_{i=1}^{k} \frac{\log_2 x}{p_i} + O\left(k \frac{(\log_3 x)^2}{\log_2 x}\right) = \sum_{i=1}^{k} \frac{\log_2 x}{p_i} + O\left(\frac{(\log_3 x)^3}{\log_2 x}\right).\end{align*}
Therefore,
\begin{align*} \prod_{\substack{z < q\le x^{1/\log_2 x} \\ q\equiv 1\!\!\!\!\!\pmod{p_i}\,\text{for some $i$}}}\left(1-\frac{1}{q}\right) &= \left(\prod_{i=1}^{k} \exp\left(-\frac{\log_2 x}{p_i}\right)\right) \Bigg(1+O\bigg(\frac{(\log_3 x)^4}{\log_2 x}\bigg)\Bigg) \\
&= \prod_{i=1}^{k}\exp\left(-\frac{\log_2 x}{p_i}\right)  + O\bigg(\frac{(\log_3 x)^4}{\log_2 x}\bigg). \end{align*}

Now collect estimates. We find that the number of $n \in \Ss_k$ where $p_1,\dots,p_k$ are the prime divisors of $n$ from $(y,z]$ is
\begin{equation}\label{eq:sktosum} x \frac{e^{-\gamma}}{\log{z}} \left(\frac{1}{p_1\cdots p_k}-\prod_{i=1}^{k} \frac{\exp(-(\log_2 x)/p_i)}{p_i}\right) + O\left(\frac{x}{p_1\cdots p_k \exp(K_5\sqrt{\log_3 x})} \right). \end{equation}

Finally, we sum \eqref{eq:sktosum} over all sets of distinct primes $p_1,\dots,p_k \in (y,z]$. The $O$-terms contribute $O(x/\exp(K_5\sqrt{\log_3 x}))$. Next we look at the contribution from the $1/p_1\cdots p_k$ terms. On the one hand, the multinomial theorem immediately implies that
\[ \sum_{y < p_1 < p_2 < \dots < p_k \le z} \frac{1}{p_1\cdots p_k} \le \frac{1}{k!}\sigma_0^k,\quad\text{where}\quad \sigma_0 := \sum_{y < p \le z}\frac{1}{p}. \]
(We have $\sigma_0 \asymp 1/\sqrt{\log_3 x}$ for large $x$ by Mertens' theorem.) On the other hand,
\begin{align*} \sum_{\substack{p_1,\dots,p_k \in (y,z] \\ \text{distinct} }} \frac{1}{p_1\cdots p_k} &= \sum_{\substack{p_1,\dots,p_{k-1} \in (y,z]\\ \text{distinct}}} \frac{1}{p_1\cdots p_{k-1}}\sum_{\substack{y < p_k \le z \\ p_k\notin \{p_1,\dots,p_{k-1}\}}}  \frac{1}{p_k} \\ &\ge \left(\sigma_0-\frac{k-1}{y}\right) \sum_{\substack{p_1,\dots,p_{k-1} \in (y,z] \\ \text{distinct}}} \frac{1}{p_1\cdots p_{k-1}}.
 \end{align*}
We can estimate the sum over $p_1,\dots,p_{k-1}$ in a similar way. Iterating, we find that
\[ \sum_{\substack{p_1,\dots,p_{k} \in (y,z] \\ \text{distinct}}}  \frac{1}{p_1\cdots p_k} \ge \prod_{i=0}^{k-1}\left(\sigma_0-\frac{i}{y}\right)\ge \left(\sigma_0 - \frac{2(\log_3 x)^2}{\log_2 x}\right)^k, \]
so that
\[ \sum_{y < p_1 < p_2 < \dots < p_k \le z} \frac{1}{p_1\cdots p_k} \ge \frac{1}{k!}\left(\sigma_0 - \frac{2(\log_3 x)^2}{\log_2 x}\right)^k. \]
Combining the upper and lower bounds,
\[ \sum_{y < p_1 < p_2 < \dots < p_k \le z} \frac{1}{p_1\cdots p_k} = \frac{1}{k!} \sigma_0^k \left(1 + O\left(\frac{(\log_3 x)^3}{\log_2 x}\right)\right)^k = \frac{1}{k!} \sigma_0^k + O\left(\frac{1}{k!}\frac{(\log_3 x)^4}{\log_2 x}\right).\]

The contribution from the terms of the form $\prod_{i=1}^{k} \exp(-(\log_2 x)/p_i)/p_i$ can be handled similarly. Put
\[ \sigma_1 := \sum_{y < p \le z} \frac{\exp(-(\log_2 x)/p)}{p}.  \]
Clearly, $\sigma_1 \le \sum_{y < p \le z} 1/p \ll 1/\sqrt{\log_3 x}$. Since $\exp(-(\log_2 x)/p) \gg 1$ when $p \ge \log_2 x$, we also have that $\sigma_1 \gg \sum_{\log_2 x < p \le z} 1/p \gg 1/\sqrt{\log_3 x}$. Now a computation completely parallel to the one shown above yields
\[ \sum_{y < p_1 < p_2 < \dots < p_k \le z} \prod_{i=1}^{k} \frac{\exp(-(\log_2 x)/p_i)}{p_i} = \frac{1}{k!} \sigma_1^k + O\left(\frac{1}{k!}\frac{(\log_3 x)^4}{\log_2 x}\right). \]

Piecing everything together, we conclude that
\begin{equation}\label{eq:skestimatefinal} \#\Ss_k = e^{-\gamma} \frac{x}{\log{z}}\left(\frac{\sigma_0^k}{k!} - \frac{\sigma_1^k}{k!}\right) + O\left(\frac{x}{\exp(K_5\sqrt{\log_3 x})} + \frac{x}{k!}\frac{(\log_3 x)^4}{\log_2 x}\right).\end{equation}
\subsection{Denouement}\label{sec:DayNewMa} Summing \eqref{eq:skestimatefinal} over positive integers $k\le \log_3 x$, keeping in mind that $\sigma_0, \sigma_1 \ll 1/\sqrt{\log_3 x}$, we find that
\[ \sum_{1\le k \le \log_3 x}\#\Ss_k = e^{-\gamma}\frac{x}{\log{z}} (\exp(\sigma_0)-\exp(\sigma_1)) + O\left(\frac{x}{\exp(K_6 \sqrt{\log_3 x})}\right).\]
By Mertens' theorem, $\exp(\sigma_0) = \frac{\log{z}}{\log{y}}\left(1 + O(1/\exp(K_7\sqrt{\log_3 x}))\right)$. So recalling \eqref{eq:s0estimate},
\[ \#\Ss_0 - \sum_{1\le k \le \log_3 x}\#\Ss_k = e^{-\gamma}\frac{x}{\log{z}} \exp(\sigma_1) + O(x/\exp(K_8 \sqrt{\log_3 x})). \]
By another application of the prime number theorem with the de la Vall\'ee Poussin error term,
\[ \sigma_1 = \int_{y}^{z} \frac{\exp(-(\log_2 x)/t)}{t\log{t}} \,d\theta(t) = \int_{y}^{z} \frac{\exp(-(\log_2 x)/t)}{t\log{t}} \,dt + O(1/\exp(K_9 \sqrt{\log_3 x})),\]
and thus
\begin{equation}\label{eq:sdiff} \#\Ss_0 - \sum_{1\le k \le \log_3 x}\#\Ss_k = e^{-\gamma} \frac{x}{\log{z}} \exp\left(\int_{y}^{z} \frac{\exp(-(\log_2 x)/t)}{t\log{t}} \,dt\right) + O(x/\exp(K_{10} \sqrt{\log_3 x})). \end{equation}

We proceed to analyze the integral appearing in this last estimate. Making the change of variables $u=(\log_2 x)/t$,
\[ \int_{y}^{z} \frac{\exp(-(\log_2 x)/t)}{t\log{t}} \,dt = \frac{1}{\log_3 x}\int_{(\log_2 x)/z}^{\log_3{x}} \frac{\exp(-u)}{u}\left(1-\frac{\log{u}}{\log_3 x}\right)^{-1}\, du.\]
Inside the domain of integration, $\log{u} \ll \sqrt{\log_3 x}$, and so for each fixed positive integer $M$,
\[ \left(1-\frac{\log{u}}{\log_3 x}\right)^{-1} = 1 +\left(\frac{\log{u}}{\log_3 x}\right)+\left(\frac{\log{u}}{\log_3 x}\right)^2+\dots + \left(\frac{\log{u}}{\log_3 x}\right)^M + O_M((\log_3 x)^{-(M+1)/2}).\]
Thus,
\begin{multline*} \frac{1}{\log_3 x}\int_{(\log_2 x)/z}^{\log_3{x}} \frac{\exp(-u)}{u}\left(1-\frac{\log{u}}{\log_3 x}\right)^{-1}\, du \\ = \sum_{k=0}^{M} \frac{1}{(\log_3 x)^{k+1}} \int_{(\log_2 x)/z}^{\log_3{x}} \frac{\exp(-u)}{u} \log^{k}{u}\, du \\+ O\left(\frac{1}{(\log_3 x)^{(M+3)/2}} \int_{(\log_2{x})/z}^{\log_3{x}} \frac{\exp(-u)}{u}\, du\right).
\end{multline*}
The $O$-term here is $\ll (\log_3 x)^{-\frac{1}{2}(M+3)} \int_{(\log_2 x)/z}^{\log_3 x} du/u \ll (\log_3 x)^{-1-\frac{1}{2}M}$. To handle the main term, we integrate by parts to find that
\begin{multline*} \int_{(\log_2 x)/z}^{\log_3{x}} \frac{\exp(-u)}{u} \log^{k}{u}\, du  = \left.\exp(-u)\frac{\log^{k+1}{u}}{k+1}\right|_{u=(\log_2 x)/z}^{u=\log_3 x} \\ + \frac{1}{k+1}\int_{(\log_2 x)/z}^{\log_3 x} \exp(-u)\log^{k+1}{u}\, du.\end{multline*}
For each $0 \le k \le M$, and all large $x$,
\[ \left.\exp(-u)\frac{\log^{k+1}{u}}{k+1}\right|_{u=(\log_2 x)/z}^{u=\log_3 x} \\
= \frac{-1}{k+1}\log\left(\frac{\log_2 x}{z}\right)^{k+1} + O_M(1/\exp(K_{11}\sqrt{\log_3 x})), \]
while
\begin{align*} \frac{1}{k+1}\int_{(\log_2 x)/z}^{\log_3 x} &\exp(-u)\log^{k+1}{u}\, du \\
&= \frac{1}{k+1}\int_{0}^{\infty} \exp(-u)\log^{k+1}{u} \, du + O_M(1/\exp(K_{12}\sqrt{\log_3 x})) \\
%&= \frac{1}{k+1}\left.\left(\frac{d^n}{dt^{k+1}}\int_{0}^{\infty} u^{t-1} \exp(-u)\, du\right)\right|_{s=1}+  O_M(1/\exp(K_{10}\sqrt{\log_3 x}))\\
&= \frac{1}{k+1} \Gamma^{(k+1)}(1) + O_M(1/\exp(K_{12}\sqrt{\log_3 x})) \\ &= k! C_{k+1} + O_M(1/\exp(K_{12}\sqrt{\log_3 x})). \end{align*}

Assembling our results,
\begin{align*} \int_{y}^{z}& \frac{\exp(-(\log_2 x)/t)}{t\log{t}} \,dt \\ &\hphantom{skip}= -\sum_{k=0}^{M} \frac{1}{k+1} \left(\frac{\log((\log_2 x)/z)}{\log_3 x}\right)^{k+1}
 + \sum_{k=0}^{M} \frac{k! C_{k+1}}{(\log_3 x)^{k+1}} + O_M((\log_3 x)^{-1-\frac{1}{2}M})\\
 &\hphantom{skip}= \log\left(1-\frac{\log((\log_2 x)/z)}{\log_3 x}\right)+   \sum_{k=0}^{M} \frac{k! C_{k+1}}{(\log_3 x)^{k+1}} + O_M((\log_3 x)^{-1-\frac{1}{2}M}) \\
 &\hphantom{skip}= \log \frac{\log z}{\log_{3}x} +   \sum_{k=0}^{M} \frac{k! C_{k+1}}{(\log_3 x)^{k+1}} + O_M((\log_3 x)^{-1-\frac{1}{2}M}).\end{align*}
We now choose $M=2N$, where $N$ is as in Theorem \ref{thm:main}.
In the last displayed sum on $k$, the terms of the sum with $k\ge N$ may be absorbed into the error. Doing so and exponentiating,
\begin{multline*} \exp\left(\int_{y}^{z}\frac{\exp(-(\log_2 x)/t)}{t\log{t}} \,dt\right) \\
= \frac{\log{z}}{\log_3 x} \exp\left(\sum_{1 \le k \le N} \frac{(k-1)! C_{k}}{(\log_3 x)^{k}}\right)\left(1 + O_N((\log_3 x)^{-1-N})\right),\end{multline*}
so that
\begin{align*} e^{-\gamma} \frac{x}{\log{z}} \exp\bigg(\int_{y}^{z} &\frac{\exp(-(\log_2 x)/t)}{t\log{t}} \,dt\bigg) \\
&= e^{-\gamma} \frac{x}{\log_3 x}\exp\left(\sum_{1 \le k \le N} \frac{(k-1)! C_k}{(\log_3 x)^k}\right) \left(1 + O_N((\log_3 x)^{-1-N})\right)\\
&= e^{-\gamma} \frac{x}{\log_3 x}\exp\left(\sum_{1 \le k \le N} \frac{(k-1)! C_k}{(\log_3 x)^k}\right) + O_N(x(\log_3 x)^{-2-N}).
\end{align*}
This expression describes $\#(\Ss_0\setminus \bigcup_{1 \le k \le \log_3 x} \Ss_k)$, by \eqref{eq:sdiff}, and so also describes $C(x)$, by the discussion in \S\ref{sec:outline}.  Theorem \ref{thm:main} follows, along with the description of the constants $c_k$ appearing in the introduction.

\providecommand{\bysame}{\leavevmode\hbox to3em{\hrulefill}\thinspace}
\providecommand{\MR}{\relax\ifhmode\unskip\space\fi MR }
\providecommand{\MRhref}[2]{%
  \href{http://www.ams.org/mathscinet-getitem?mr=#1}{#2}
}
\providecommand{\href}[2]{#2}

\end{document}